\documentclass[conference]{IEEEtran}
\IEEEoverridecommandlockouts
\usepackage{cite}

\usepackage{amsmath,amssymb,amsfonts}
\usepackage{algorithm}
\usepackage{algcompatible}
\usepackage{placeins}
\usepackage{graphicx}
\usepackage{textcomp}
\usepackage{xcolor}
\usepackage[margin=0.75in]{geometry}

\def\BibTeX{{\rm B\kern-.05em{\sc i\kern-.025em b}\kern-.08em
    T\kern-.1667em\lower.7ex\hbox{E}\kern-.125emX}}
    
\DeclareMathOperator{\Tr}{Tr}

\begin{document}

\title{Equitable dynamic electricity pricing via implicitly constrained dual
and subgradient methods\\
}

\author{\IEEEauthorblockN{Emmanuel Balogun}
\IEEEauthorblockA{\textit{Department of Mechanical Engineering} \\
\textit{Stanford University}\\
Stanford, USA \\
ebalogun@stanford.edu}
\and
\IEEEauthorblockN{Sonia Martin}
\IEEEauthorblockA{\textit{Department of Mechanical Engineering} \\
\textit{Stanford University}\\
Stanford, USA \\
soniamartin@stanford.edu}
\and
\IEEEauthorblockN{Anthony Degleris}
\IEEEauthorblockA{\textit{Department of Electrical Engineering} \\
\textit{Stanford University}\\
Stanford, USA \\
degleris@stanford.edu}
}

\author{
    \IEEEauthorblockN{Emmanuel Balogun\IEEEauthorrefmark{1}, Sonia Martin\IEEEauthorrefmark{1}, Anthony Degleris\IEEEauthorrefmark{2}, Ram Rajagopal\IEEEauthorrefmark{3}}
    \IEEEauthorblockA{\IEEEauthorrefmark{1}Department of Mechanical Engineering, Stanford University, Stanford, USA 
    }
    
    \IEEEauthorblockA{\IEEEauthorrefmark{2}Department of Electrical Engineering, 
    Stanford University, Stanford, USA}
    \IEEEauthorblockA{\IEEEauthorrefmark{3}Department of Civil and Environmental Engineering, 
    Stanford University, Stanford, USA
    }
}

\maketitle

\begin{abstract}
 Coordination of distributed energy resources is critical for electricity grid management. Although nodal pricing schemes can mitigate congestion and voltage deviations, the resulting prices are not necessarily equitable. In this work, we leverage market mechanisms for DER coordination and propose a daily dynamic nodal pricing scheme that incorporates equity. We introduce a pricing ``oracle,'' which we call the Power Distribution Authority, that sets equitable prices to manage the grid. We present two algorithms for executing this scheme and show that both methods are able to set prices that satisfy both voltage and equity constraints. Both proposed algorithms also outperform the common utility time-of-use pricing schemes by at least 45\%. New market mechanisms are needed as the grid is transforming, and power system operators may leverage these methods for pricing electricity in a grid-aware, equitable fashion.
\end{abstract}

\begin{IEEEkeywords}
Dynamic pricing, Equity, Electric Vehicles, Optimization
\end{IEEEkeywords}


\section{Introduction}
The electricity grid is experiencing a tectonic shift toward decentralized power generation, driven by the rapid adoption of distributed energy resources (DERs) such as solar photovoltaics, electric vehicles (EVs), and batteries. The accelerating rate of EV adoption is leading to a large growth in peak power demand. Deploying new distribution or transmission infrastructure to accommodate peak demands is usually time-consuming (due to planning and permits) and expensive, so coordinating DERs to reduce peak demand is essential. A lack of EV charging coordination can have significant impacts on electricity infrastructure, even when EV adoption is low~\cite{muratori2018impact}. DER coordination within a distribution grid is challenging because DERs are often owned by separate, profit-driven entities who act without regard to grid impacts, often because they are not incentivized to do so. In this paper we propose methods to indirectly coordinate DERs by influencing electricity demand via electricity prices.

In markets where participants are economically rational and sufficiently flexible, dynamic pricing can shape demand by effectively coordinating resources using prices~\cite{denboer_2015, luo2017stochastic, kim2014dynamic}. Dynamic nodal pricing can enact demand profile changes on the nodal (or grid bus) level by varying nodal prices over time horizons~\cite{andrianesis2020distribution}. For example, customers can be discouraged from charging EVs at times of congestion due to high prices, thus reducing peak demand, and demand response programs can be controlled via dynamic pricing~\cite{nguyen2016dynamic}. However, dynamic nodal pricing can be unfair. Customers connected to nodes with lower hosting capacities and fewer behind-the-meter (BTM) resources may end up paying more for energy than customers at favorable nodes, further exacerbating inequity~\cite{brockway2021inequitable}. Furthermore, one study finds that in a case of dynamic residential pricing, assuming that customers did not change their usage patterns, only 35\% of customers saw cheaper bills, showing that there is a need to set dynamic prices in an explicitly equitable manner~\cite{horowitz2014equity}.

Many prior works have proposed algorithms to derive dynamic prices, as well as improve their stability. Dynamic locational marginal prices (LMPs), or grid nodal prices, can be calculated via various descent methods, such as iterative gradient descent~\cite{okawa2014dynamic} and subgradient descent~\cite{roozbehani2010dynamic}. Incorporating extreme price reduction is one method to stabilize the LMPs in a system~\cite{roozbehani2010dynamic}. Other studies examine the role of energy storage systems and dynamic pricing~\cite{kang2021distributed}, as well as EVs and dynamic pricing~\cite{luo2017stochastic, moghaddam2019coordinated}. EVs can be leveraged and coordinated via a dynamic pricing scheme to help balance loads within a power grid~\cite{moghaddam2019coordinated}. To further improve coordination, an EV service provider that owns a set of EV chargers can profit from setting dynamic prices to smooth out uncertainties in consumer charging behavior~\cite{luo2017stochastic}. 


In addition to coordinating EV charging to balance loads, dynamic pricing can also help mitigate voltage fluctuations, which is important to protect grid-connected devices within the distribution grid~\cite{chattopadhyay2001pricing}. While uncoordinated EV charging can lead to voltage drops below the regulation minimum~\cite{deilami2011real}, controlling EV charging with dynamic pricing can reduce voltage violations caused by other DERs~\cite{then2023coordinated}. The power distribution grid's sensitivity to undesired voltage deviations necessitates attentive design of dynamic pricing: a pricing scheme should take into account voltage stability and impact on other grid assets.


While prior works propose separate solutions for dynamic pricing and reducing load-induced voltage violations, they do not explicitly prioritize voltage deviations and economic/price equity within a unifying framework. This work proposes and evaluates dynamic pricing methods to minimize voltage violations while preserving price equity. Specifically, we leverage the dual and the subgradient descent methods. The rest of the paper is organized as follows: in Section~\ref{methodology}, we describe the methods and algorithms, including problem setup; in Section~\ref{results}, we present and discuss the simulation results for the proposed algorithms, and finally in Section~\ref{conclusion}, we conclude and discuss future work.

\section{Methodology}\label{methodology}
              Excessive congestion or overloading in certain portions of the power grid may lead to severe voltage fluctuations and violations~\cite{limmer2019peak, moghaddam2019coordinated}. The methods in this work aim to reduce voltage violations/fluctuations while preserving price equity. In this section, we first describe the data used for the simulation and then provide an overview of the power flow equations. Next, we introduce the Power Distribution Authority (PDA) ``oracle,'' which simulates charging behavior in order to set electricity prices. We then define price equity and voltage fluctuations within the grid. Lastly, we describe the proposed pricing algorithms.

\subsection{Problem Setup and Data}\label{sec:setup}
To demonstrate the price-setting algorithm, we analyze a 14-bus test case with generator and load buses. A set of EV charging stations are located at some of the load buses, along with collocated stationary batteries. Given a load forecast of the EV chargers for the next day, the output of the algorithm is a set of electricity prices for each node.

\begin{figure}[t]
    \centering
    \includegraphics[width=.8\linewidth]{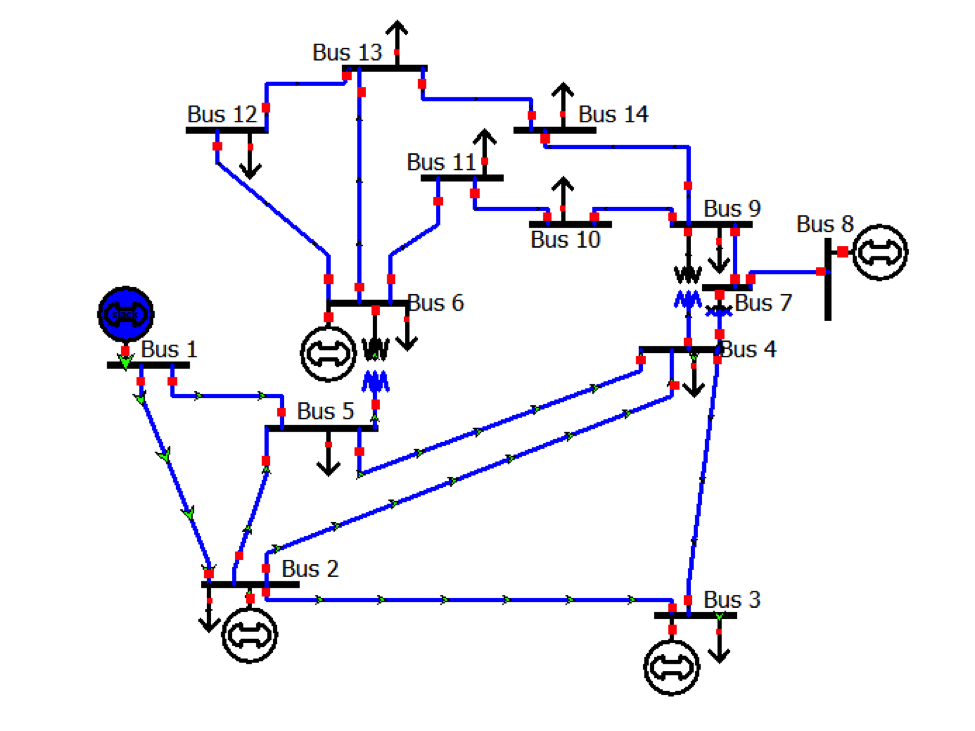}
    \caption{Diagram of the IEEE 14-bus test feeder network~\cite{bus_data}. EV charging stations are located on buses 11-14.}
    \label{fig:test_feeder}
\end{figure}

 Two datasets provide the foundation of the simulation: grid system information and EV charging profile data. We use the IEEE 14-bus test case, shown in Fig.~\ref{fig:test_feeder}, which portrays a synthetic power network with 14 nodes~\cite{bus_data}. This provides network information, including the line admittance data, which are inputs to the power flow equations. It also gives locations and magnitudes of spot loads and generators, which are existing loads and power sources within the grid under normal loading conditions. Bus 1 is the slack bus. We site EV charging stations at 4 of the 14 buses of the power network. For each station, we leverage the SPEECh model~\cite{powell2022scalable} to generate time-varying EV charging loads based on a population of 3200 EVs.

\subsection{Power Flow Problem}
In a power system, the real and reactive powers at every node are defined, respectively, as follows~\cite{glover2012power}:
\begin{align}
    P_i &= \sum_{k = 1}^{n}|v_i||v_k|\left[G_{ik}\cos{(\theta_{i}-\theta_k)} + B_{ik}\sin{(\theta_{i}-\theta_k)}\right]\label{P_pf} \\
    Q_i &= \sum_{k = 1}^{n}|v_i||v_k|\left[G_{ik}\sin(\theta_{i}-\theta_k) - B_{ik}\cos(\theta_{i}-\theta_k)\right]\label{Q_pf}
\end{align}
where $n$ is the number of nodes within the power network, $\theta_i - \theta_k$ is the voltage phase angle difference between bus $i$ and bus $k$, $G_{ik}$ is the line conductance between the buses, $B_{ik}$ is the line susceptance, and $P_i$ and $Q_i$ are the net real and reactive power, respectively, generated at node $i$.
 Using the Newton-Raphson method~\cite{glover2012power}, we can linearize the power flow equations to obtain a Taylor series approximation about a stationary operating point, similar to the process in~\cite{luo2017stochastic}. This process makes the subtle assumptions that the deviations are small enough that the approximation is useful. The initial stationary operating point is the nominal voltage magnitudes and phase angles of the test case, provided in the network information, while each subsequent stationary point is the prior time's operating point. To perform the linearization, we take the Jacobian of the system, which can be broken down into four block matrices (see Equation \eqref{Jacob}), by taking the partial derivative of both Equations~\eqref{P_pf} and~\eqref{Q_pf} with respect to voltage magnitude $v$ and phase angle $\theta$. 
 
\begin{equation}\label{Jacob}
\mathbf{J}=\begin{bmatrix}
    \frac{\partial{\mathbf{P}}}{\partial{\mathbf{V}}} & \frac{\partial{\mathbf{P}}}{\partial{\mathbf{\theta}}}\vspace{.8em} \\ 
    \frac{\partial{\mathbf{Q}}}{\partial{\mathbf{V}}} & \frac{\partial{\mathbf{Q}}}{\partial{\mathbf{\theta}}}
\end{bmatrix}
\end{equation}

Each block of the Jacobian is defined below:
 \begin{align}
     \frac{\partial{P_i}}{\partial{v_j}}& = \begin{cases}
     2v_iG_{ii} + \sum_{k=1; k \neq
     i}^{n}v_k[G_{ik}\cos{(\theta_{i}-\theta_k)} + 
     \\ \hspace{1.5em} B_{ik}\sin{(\theta_{i}-\theta_k)}]; \ j = i \\
     v_i\left[G_{ij}\cos{(\theta_{i}-\theta_j)} + B_{ij}\sin{(\theta_{i}-\theta_j)}\right];\\ \hspace{1.5em}  \ j \neq i
		 \end{cases} \\
     \frac{\partial{Q_i}}{\partial{v_j}} &= \begin{cases}
     -2v_iB_{ii} + \sum_{k=1; k\neq i}^{n}v_k[G_{ik}\sin{(\theta_{i}-\theta_k)}  \\ \hspace{1.5em} -B_{ik}\cos{(\theta_{i}-\theta_k)}]; \  j = i \\
     v_i\left[G_{ij}\sin{(\theta_{i}-\theta_j)} -B_{ij}\cos{(\theta_{i}-\theta_j)}\right]; \\\hspace{1.5em} j \neq i
     \end{cases} \\
    \frac{\partial{P_i}}{\partial{\theta_j}} &= \begin{cases}
    v_i \sum^n_{k=1;k\neq i}v_k [-G_{ik}\sin(\theta_{i}-\theta_k) \\
    \hspace{1.5em} +B_{ik}\cos(\theta_{i}-\theta_k)]; \ j=i \\
    v_i v_j [G_{ij}(\sin(\theta_{i}-\theta_j)-B_{ij}\cos(\theta_{i}-\theta_j)]; \\  \hspace{1.5em} j\neq i
    \end{cases} \\
     \frac{\partial{Q_i}}{\partial{\theta_j}} &=\begin{cases}
     v_i \sum^n_{k=1;k\neq i}v_k [G_{ik}(\cos(\theta_i-\theta_k)  \\ \hspace{1.5em}
     +B_{ik}\sin(\theta_i-\theta_k)];\ j=i \\ 
    v_i v_j [-G_{ij}\cos(\theta_i-\theta_j)  \\ \hspace{1.5em}
    -B_{ij}\sin(\theta_i-\theta_j)]; \ j\neq i 
    \end{cases}
\end{align}

The first order Taylor series approximation about a fixed operating point is
\begin{equation}\label{taylor1}
    \begin{bmatrix}
        \mathbf{P} \\ \mathbf{Q}
    \end{bmatrix} - 
    \begin{bmatrix}
        \mathbf{P_{0}} \\ 
        \mathbf{Q_{0}}
    \end{bmatrix} = \mathbf{J} \left(
    \begin{bmatrix}
        \mathbf{V} \\ 
        \mathbf{\theta}
    \end{bmatrix} -
    \begin{bmatrix}
        \mathbf{V_{0}} \\ 
        \mathbf{\theta_{0}}
    \end{bmatrix} \right),
\end{equation}
where $\mathbf{P}$, $\mathbf{Q}$ are the powers from the first timestep to the last timestep and $\mathbf{P}_0$, $\mathbf{Q}_0$ are powers from the prior timesteps. This can be simplified as 
\begin{align}
    \begin{bmatrix} \Delta{\mathbf{P}}  \\ \Delta{\mathbf{Q}}\end{bmatrix} =  \mathbf{J} \begin{bmatrix} \Delta{\mathbf{V}}  \\ \Delta{\mathbf{\theta}}\end{bmatrix} .\label{eq:taylor}
\end{align}
We delete the row $\mathbf{P}_{\theta, slack}$ corresponding to the slack bus angle $\Delta \theta_{slack}$, as this is usually fixed at $0$, and redefine $\mathbf{J}$ as the Jacobian without row $\mathbf{P}_{\theta, slack}$. From Equation~\eqref{eq:taylor}, we can calculate the ``voltage deviations,'' which are defined as $[\Delta{\mathbf{V}} \ \Delta{\mathbf{\theta}}]^{T}$ by pseudo-inverting the Jacobian. We desire to minimize these deviations, expressed as 
\begin{align} 
 \left\| \mathbf{J}^{\dag} \begin{bmatrix} \Delta{\mathbf{P}}  \\ \Delta{\mathbf{Q}}\end{bmatrix} \right\|_2^2.
\end{align}

The vector $[\Delta{\mathbf{P}}\ \Delta{\mathbf{Q}}]^T$ represents change between the current and prior net load at each node.

\subsection{Power Distribution Authority}

\begin{figure}[htbp]
    \centering
    \includegraphics[width=.9\linewidth]{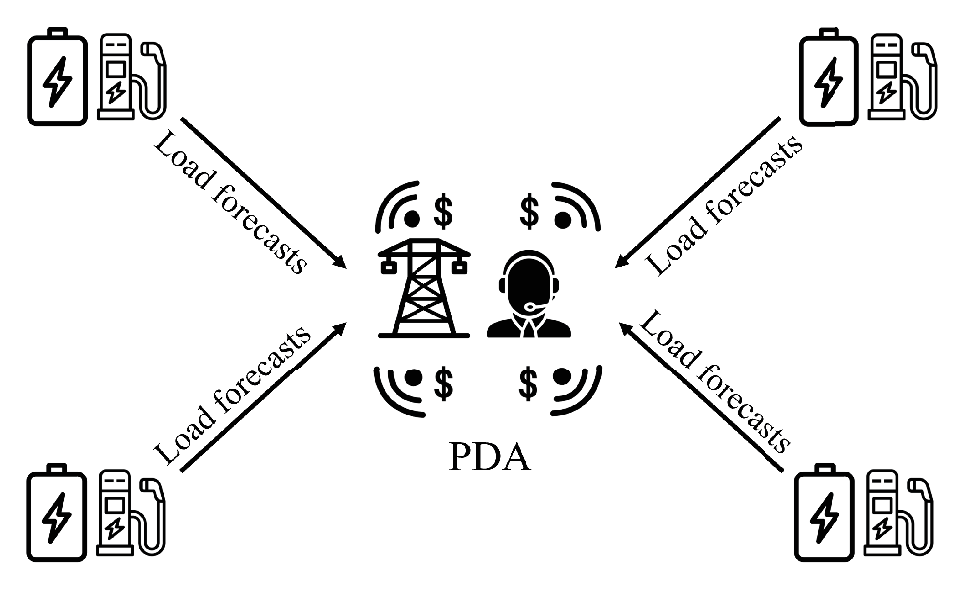}
    \caption{Power Distribution Authority diagram. EV charging station providers send their load forecasts to the PDA, who communicates dynamic nodal prices back to the chargers.}
    \label{fig:PDA}
\end{figure}
\begin{figure}
    \centering
    \includegraphics[width=.75\linewidth]{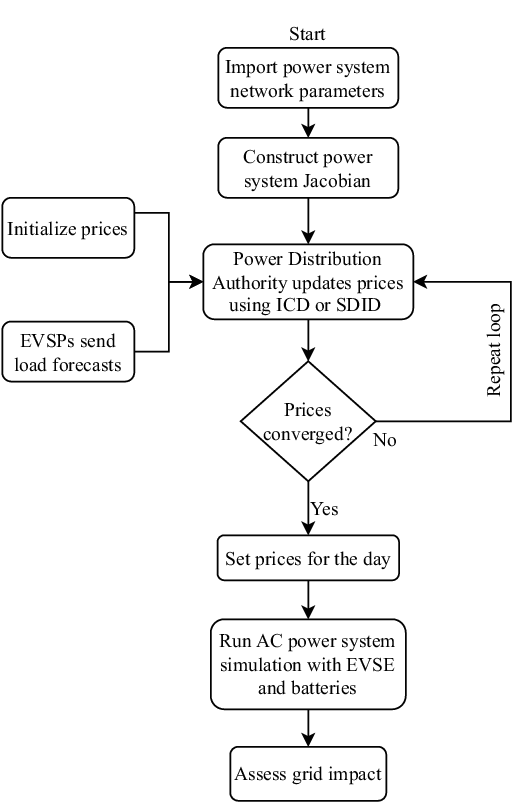}
    \caption{Algorithm for dynamic pricing scheme.}
    \label{fig:flowchart}
\end{figure}

We introduce a pricing ``oracle'' (see Fig. \ref{fig:PDA}), the Power Distribution Authority (PDA), that sets the price of electricity at each charging station node for every timestep. First, each EV charging service provider (EVSP) sends the PDA its charging forecasts and initial state-of-charge (SoC). Next, as shown in Fig.~\ref{fig:flowchart}, the PDA simulates how a profit maximizing EVSP would operate their stationary battery given the current set of time-varying prices at its node. The PDA then leverages the proposed algorithms in this paper to set the prices for the next 24 hours, fulfilling the objectives of minimal voltage violations and price equity. Using the new set of prices, this process is repeated until there is price convergence or no further improvement in the objective.



\subsection{Equity}
In this section, we mathematically express the concept of price equity. We define economic burden below as:
\begin{equation}
\label{equ1}
    \Lambda = \sum_{t=1}^{T}\lambda(t),
\end{equation}
which is the price an actor pays per kWh of energy consumed for the entire simulated timeframe and where $\lambda$ is the electricity price (\$/kWh) vector. In the context of the proposed dynamic pricing, the goal is to minimize the variance of the expected economic burden of the actors/agents within the power network, so we propose the following. The penalty controlling equity will be defined as:
\begin{align}
\Gamma(\lambda) = \underset {k \in K} {\min}  \left\{{\sum_{t=1}^{T}\lambda_{k}(t)}\right\} -  \underset {k \in K} {\max}\left\{{\sum_{t=1}^{T}\lambda_{k}(t)}\right\},
\end{align}

where $K$ is the set of charging stations and $\lambda$ is the \$/kWh paid by each actor over the horizon. Equity is bounded by $(0, -\infty)$, with higher equity level monotonically increasing with its value. It is worth noting that any convex definition of equity will be valid. For example, with regional socioeconomic data, one can refine the definition to be in terms of energy burden, which is the percentage of income spent by a customer on their energy bill \cite{hernandez2010energy}. If desired, different variants of convex equity definitions can be leveraged.

\subsection{Implicitly Constrained Dual}
The Implicitly Constrained Dual (ICD) algorithm leverages decomposition methods. We decompose the problem into two main objectives and include a consensus constraint:

\begin{align}
& \underset{\mathbf{P}, \mathbf{\hat{P}}, \mathbf{Q}}{\text{minimize}}
& & \left\| \mathbf{J}^{\dag} \begin{bmatrix} \Delta{\mathbf{P}}  \\ \Delta{\mathbf{Q}}\end{bmatrix} \right\|_2^2 \label{obj} + \Tr{\left(\mathbf{\hat{P}}_{K}^{T} \lambda_{0}\right)}  \\
& \text{subject to}
& & \Delta \mathbf{P}_{K} = \Delta \mathbf{\hat{P}}_{K}, \label{consensus}\\
& & & \mathbf{\hat{P}} = \mathbf{P}_b + \mathbf{\hat{P}}_{K}.
\end{align}

The objectives are the first and second terms included in Equation~\eqref{obj}. The first term in Equation~\eqref{obj} represents the voltage fluctuations and the second term is the EVSP electricity cost. $\Tr$ is the trace, $\Delta{\mathbf{P}} = \mathbf{P} - \mathbf{P}_{0}$, and $\Delta{\mathbf{Q}} = \mathbf{Q} - \mathbf{Q}_{0}$. Equation~\eqref{consensus} is the consensus constraint, which ensures that $\Delta \mathbf{P}_{K}$, a subset of the rows in $\Delta \mathbf{P}$, is equal to $\Delta \hat{\mathbf{P}}_{K}$, the change in time-varying load at charging buses represented by set $K$. $\mathbf{P}_b$ and $\mathbf{P}_{K}$ are the battery power and total load forecast, respectively. $\lambda_0$ is the initial electricity price matrix made up of each charging node's price. The electricity prices at each EV charging node correspond to the dual variables of the consensus constraint in Equation~\eqref{consensus}. To satisfy the equity objective, we analytically resolve the optimal dual variable $\lambda^{*}$ in terms of the primal variables.

Because $\Delta \mathbf{Q} = \mathbf{0}$ (no time-varying reactive power), we can simplify Equation \eqref{obj} by retaining only $\Delta\mathbf{P}$. We define $\tilde{\mathbf{J}}$ as the submatrix of the Jacobian $\mathbf{J}$ that represents the desired columns and rewrite the objective as:
\begin{align}    
& \underset{\mathbf{\Delta P}_K}{\text{minimize}}
& & \left\| \mathbf{\tilde{J}}^{\dag} \Delta{\mathbf{P}} \right\|_2^2 \label{obj_2} + \Tr{\left(\left({\Delta \hat{\mathbf{P}}_K} + \hat{\mathbf{P}}_{0, K}\right) ^{T}\lambda_{0}\right)}  \\
& \text{subject to}
& & \Delta \mathbf{P}_{K} = \Delta \mathbf{\hat{P}}_{K}\label{consensus_2},\\
& & & \mathbf{\hat{P}}_{K} = \mathbf{P}_b + \mathbf{{P}}_{ev, K}.
\end{align}

$\hat{\mathbf{P}}_K = \Delta \hat{\mathbf{P}}_K + \hat{\mathbf{P}}_{0, K}$. We include the consensus constraint in the objective Equation~\eqref{obj_2} via a Lagrangian variable, $\lambda$, take the derivative of the augmented objective w.r.t $\Delta \mathbf{P}_K$, and resolve the optimal dual variable $\lambda^{*}$ as a function of the primal variable, expressed in Equation \eqref{equ:dual_primal_resolve}:
\begin{align}
    \label{equ:dual_primal_resolve}
    \lambda^{*}_{K}(\mathbf{P}) = \left[2(\mathbf{\tilde{J}}^{\dag})^{T}\mathbf{\tilde{J}}^{\dag}\Delta{\mathbf{P}}\right]_{K},
\end{align}

 We pass $\lambda_{K}^{*}(\mathbf{P})$ as argument to the equity function $\Gamma$, yielding the optimization problem: 
\begin{align}
& \underset{\mathbf{P}, \mathbf{\hat{P}}, \mathbf{Q}}{\text{minimize}}
& & \alpha \left\| \mathbf{J}^{\dag} \begin{bmatrix} \Delta{\mathbf{P}}  \\ \Delta{\mathbf{Q}}\end{bmatrix} \right\|_2^2 + \Tr{\left(\mathbf{\hat{P}}_{K}^{T} \lambda_{0}\right)}  - \beta \Gamma \left(\lambda_{K} ^ *(\mathbf{P})\right)\label{obj_final_icd}
\\
& \text{subject to}
& & \Delta \mathbf{P}_{K} = \Delta \hat{\mathbf{P}}_{K}, \\
& & & \mathbf{\hat{P}}_{K} = \mathbf{P}_b + \mathbf{{P}}_{ev, K},\\
& & & \mathbf{P}_t ^T \mathbf{1}= \mathbf{P}_{gen,t}^T \mathbf{1},  \hspace{1em} \label{p_gen}\\
& & & \mathbf{Q}_t^T \mathbf{1}  = \mathbf{Q}_{gen,t}^T \mathbf{1}, \hspace{1em}\label{q_gen} \\
& & &    | \mathbf{P}_{b,k} | \leq P_{b, max},\hspace{1em}  \forall k \in K,\label{batt_max_power} \\
 & & &   \mathrm{SoC}_k(t+\Delta t) = \mathrm{SoC}_k(t)+ P_{b, k}(t) \frac{\Delta t}{\mathrm{E}_k},\label{soc_evol}\notag\\
 & & & \quad \quad\forall k \in K,  \\ 
& & &    \mathrm{SoC}_{k}(t=0) =\mathrm{SoC}_{k,init} \hspace{1em},  \forall k \in K, \label{soc_init_constr}\\
& & &    0\leq \mathrm{SoC}_k(t) \leq 1 \hspace{1em},  \forall k \in K,\label{soc_range_constr}
\end{align}
where $\alpha$ and $\beta$ are weighting constants, $t$ is time, $P_{b,k}$ is the power of battery $k$ in MW, $\Delta t$ 
is the timestep length (15 minutes), $\mathrm{E}_k$ is the battery energy capacity of battery $k$ in MWh, and $\mathrm{SoC}_{k, init}$ is the initial SoC of the $k$'th battery. Equation \eqref{p_gen} and~\eqref{q_gen} are the power balance constraints, ensuring that power supply matches demand. Equation \eqref{batt_max_power} constrains the battery power, and Equation \eqref{soc_evol}--\eqref{soc_range_constr} are the battery SoC constraints.

We start the first iteration of the algorithm by randomly initializing prices, and  solving the problem. The subsequent iterations are executed with the prices $\lambda^{*}$ obtained from the dual variable of the consensus constraint of the prior iterations.
We repeat this process using the prior iteration's prices until convergence or no further objective improvement.

\subsection{Subgradient descent via implicit differentiation}
In this section, we introduce a subgradient descent algorithm that leverages implicit differentiation (SDID).
The algorithm first considers the optimization problem solved by each EVSP $k$ after prices $\lambda$ are set,
\begin{align} \label{eq:evsp-problem}
 &   \underset{\mathbf{P_k}}{\textrm{minimize}} & & 
     \mathbf{P}_k^T \mathbf{\lambda} \\
&     \textrm{subject to} & &
   \mathbf{P}_k = \mathbf{P}_{b, k} + \mathbf{P}_{ev, k}, \\
 & & &   \mathrm{SoC}_k(t+\Delta t) = \mathrm{SoC}_k(t)+ P_{b, k}(t) \frac{\Delta t}{\mathrm{E}_k}, \notag \\
& & &    | \mathbf{P}_{b, k} | \leq P_{b, max}, \hspace{1em}  \\
& & &    \mathrm{SoC}_k(t = 0) =\mathrm{SoC}_{k,init}, \hspace{1em}  \\
& & &    0\leq \mathrm{SoC}_k(t) \leq 1.
\end{align}
We now define $\mathbf P_k^*(\lambda) \in \mathbb{R}^T$ to be the optimal real power injections set by EVSP $k$ in~\eqref{eq:evsp-problem}.
We assume a power factor of 1, thus, each EVSP sets $\mathbf{Q}_k = 0$.
As before, the PDA seeks to find prices that minimize a weighted combination of voltage deviations and inequity,
\begin{align}
   F(\lambda) = 
   f(\mathbf P^*(\lambda), \lambda) =  \alpha 
   \left\| \mathbf{J}^{\dag} \begin{bmatrix} \Delta{\mathbf{P}}  \\ \Delta{\mathbf{Q}}\end{bmatrix} \right\|_2^2 - \beta \Gamma(\lambda)
\end{align}
To achieve this, the SDID algorithm iteratively updates prices using gradient steps $\lambda^{j+1} = \lambda^j - \eta \nabla F(\lambda)$, where $\eta > 0$ is a fixed step size.
Using the chain rule, we know the gradient $\nabla F(\lambda) \in \mathbb{R}^{nT}$ is,
\begin{equation}
\label{implicit_subgrad}
    \nabla F(\lambda) = \partial_1 f(\mathbf P^*(\lambda), \lambda)^T \nabla \mathbf P^*(\lambda)  + \partial_2 f(\mathbf P^*(\lambda), \lambda),
\end{equation}
where $\partial_i f(\mathbf P^*(\lambda), \lambda) \in \mathbb{R}^{nT}$ is the partial gradient of $f$ with respect to its $i$th argument
(evaluated at $(\mathbf P^*(\lambda), \lambda)$), 
and $\partial \mathbf P^*(\lambda) \in \mathbb{R}^{nT \times nT}$ is the Jacobian of $\mathbf P^*(\lambda)$, where $n$ is number of nodes and $T$ is number of timesteps. 
The function $\mathbf P^*(\lambda)$ is not explicitly defined via an analytical formula;
rather, they are implicitly defined by the solution to~\eqref{eq:evsp-problem}.
We address this by applying the implicit function theorem to a fixed-point operator whose solution corresponds to the solution of the optimization problem in~\eqref{eq:evsp-problem}, as described in~\cite{agrawal2019differentiating}.
To summarize, the SDID algorithm applies the following steps:

\begin{itemize}
    \item Randomly initialize electricity TOU prices for all EVSPs at their respective nodes
    \item Solve each EVSP profit maximization problem, with its respective constraints
    \item Obtain the gradient $\nabla F(\lambda)$ of each EVSP power w.r.t to prices
    \item Compute the subgradient of the overall objective w.r.t power and prices using Equation \eqref{implicit_subgrad}
    \item Use the implicit differentiation relationship to obtain subgradient of overall objective w.r.t prices $\lambda$ and run the update Algorithm \ref{price_update_pda}
\end{itemize}


\begin{algorithm}
\caption{SDID Price updates}\label{diffPDA}
\begin{algorithmic}
\State $N \gets \text{no. iterations}$
\State \text{decay: \{True, False\}}
\State $\gamma \gets \gamma_o$
\State $\eta \gets \eta_{init}$
\State $\lambda \gets \text{random initialization}$
\FOR{$i = 1, ..., N$}
\State $\Delta \lambda = \nabla F(\lambda) $ \COMMENT{This is done in the PDA object}
\State $\lambda \gets \lambda - \eta \Delta \lambda$
\IF{\text{decay}}
\State $\eta \gets \eta \times \gamma$
\ENDIF
\ENDFOR
\end{algorithmic}
\label{price_update_pda}
\end{algorithm}

\subsection{Power Flow Simulation}
The original power flow physics is non-convex, however, we make convex approximations to obtain prices. To ensure the voltage deviations from the problem satisfy physical power flow properties, we run a power flow simulation using an open-source 3-phase power flow solver, pandapower~\cite{pandapower.2018}. 

\section{Results and Discussion}\label{results}
\begin{figure}[h]
    \centering
    \includegraphics[width=.9\linewidth]{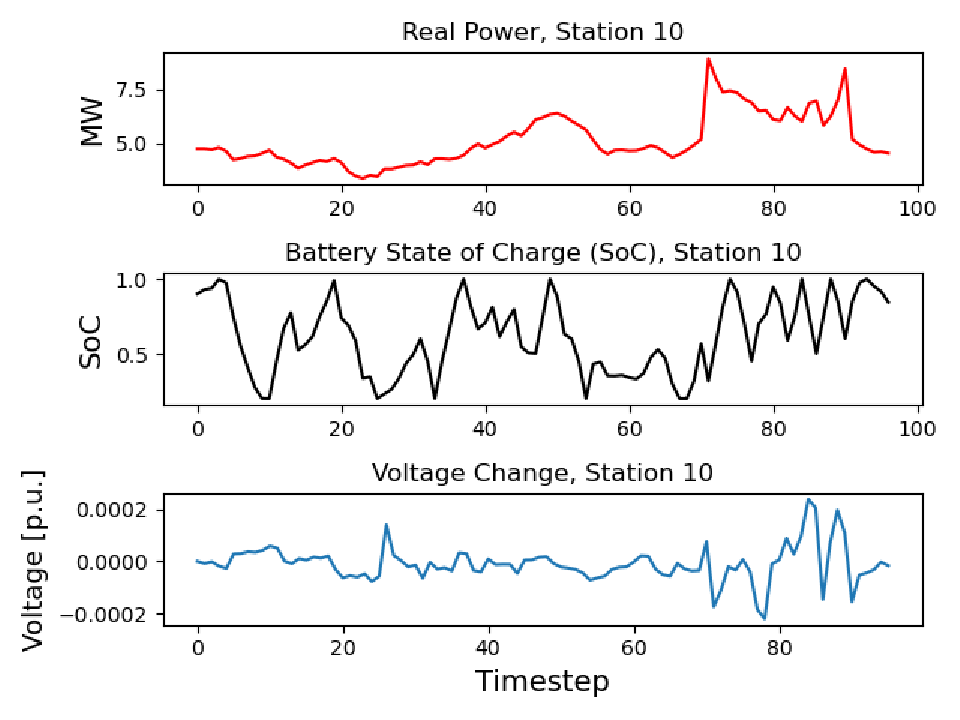}
    \caption{Plot of net power, stationary battery SoC, and voltage fluctuations for Station 10 (which is located at bus 11). The battery is able to smooth out the power profile, reducing voltage fluctuations.}
    \label{fig:station_10_ICD}
\end{figure}

We present relevant power output and convergence plots for the proposed algorithms, ICD and SDID.
The battery storage capacity at each bus can significantly impact the voltage deviations within the grid; therefore, we run two problem configurations with different battery capacities, described in Table~\ref{tab:config_1}. In both configurations, the battery at each bus has a maximum power output of 1C. We compare the prices set by the PDA to a baseline set of prices: a time-of-use (TOU) pricing scheme from PG\&E-BEV2~\cite{PGE_BEV_SCH}. Figure \ref{fig:station_10_ICD} shows the real power at bus 11 (EVSP 10), battery SoC evolution, and the approximated voltage deviations.

\begin{table}[h]
    \caption{Problem Configurations}
    \centering
    \begin{tabular}{|c|c|c|c|c|c|}
    \hline
       Config. & Variable & Bus 11 & Bus 12 & Bus 13 & Bus 14 \\
        \hline
         1 & Capacity [MWh] & 0.5 & 0.1 & 1 & 10\\
         \hline
         2 & Capacity [MWh] & 10 & 0.1 & 0.1 & 0.2\\
         \hline
        1 \& 2 &  Initial SoC & 0.9 & 0.6 & 0.6 & 0.8\\
         \hline
    \end{tabular}

    \label{tab:config_1}
\end{table}


\subsection{Implicitly Constrained Dual Method Results}
\begin{figure}[h!]
    \centering
    \includegraphics[width=.85\linewidth]{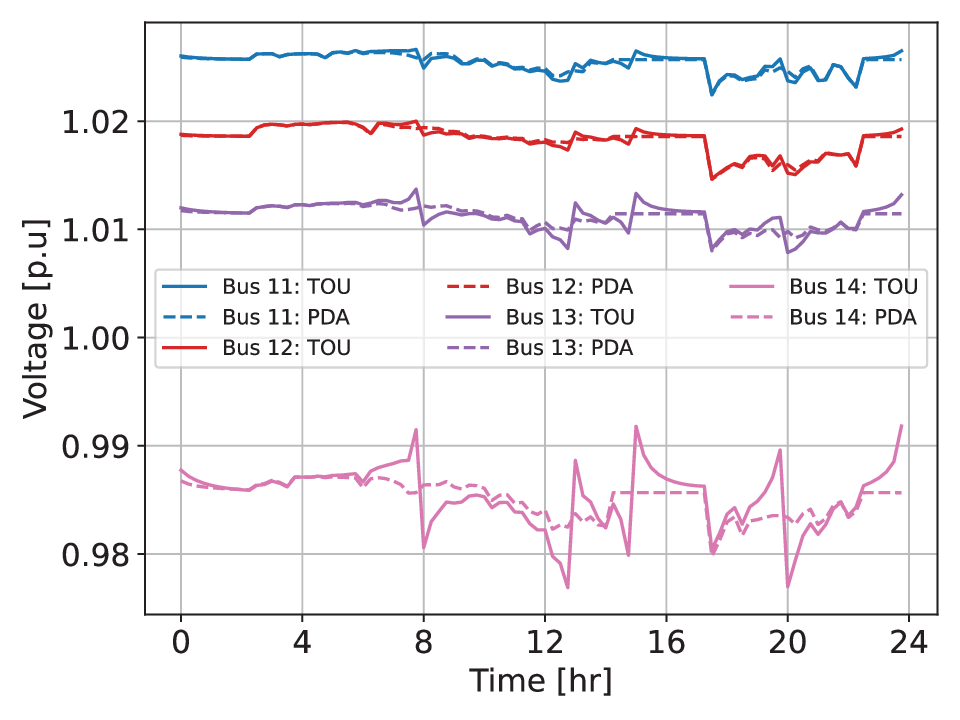}
    \caption{Voltages for buses 11-14 for ICD algorithm and configuration 1. At bus 14, the large stationary battery smooths out many of the voltage spikes relative to the TOU baseline case.}
    \label{fig:power_sim_dual_1}
\end{figure}

\begin{figure}[h!]
    \centering
    \includegraphics[width=.85\linewidth]{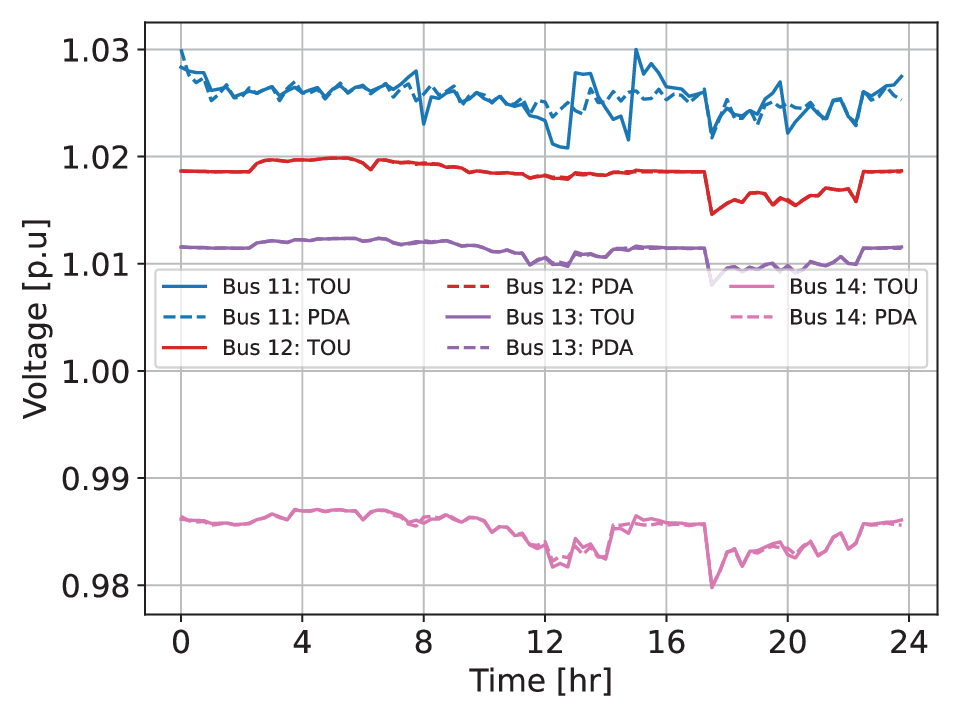}
    \vspace{-0.2cm}
    \caption{Voltages for buses 11-14 for ICD algorithm and configuration 2. At bus 14, the smaller battery is no longer able to smooth out the voltage profile, but at bus 11 the large battery has a smoothing effect.}
    \label{fig:power_sim_dual_2}
\end{figure}

The ICD algorithm is able to constrain equity to nearly 0 (as shown in Table~\ref{tab:v_compare}) and converges in one iteration.
In Figures~\ref{fig:power_sim_dual_1} and~\ref{fig:power_sim_dual_2}, the ICD method yields smoother voltage profiles than the TOU baseline. In configuration 1, bus 14 has a large battery capacity of 10 MWh, which allows it to effectively respond to prices which reduce voltage deviations. In both configurations, the prices generated by the PDA effectively reduce voltage fluctuations.

\subsection{Implicit Differentiation Method Results}
\begin{figure}[h!]
    \centering
    \includegraphics[width=.9\linewidth]{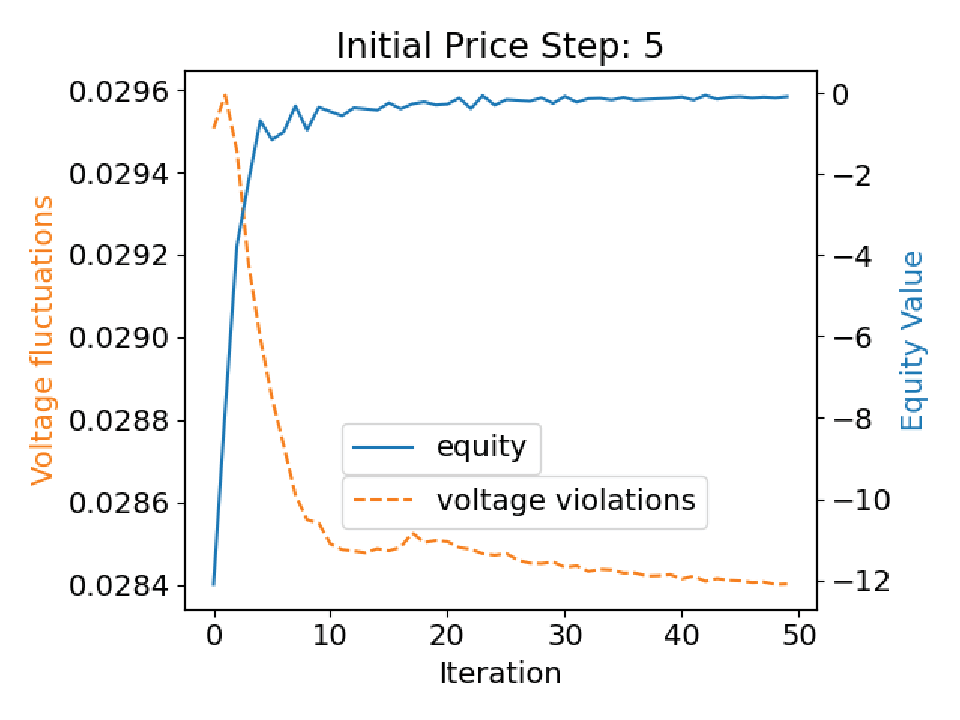}
    \caption{This plot shows equity and voltage violations for the SDID algorithm running on problem configuration 1 with a starting subgradient descent step size of 5. The algorithm simultaneously decreases voltage violations while driving equity toward zero.}
    \label{fig:ICD_training}
\end{figure}

Fig.~\ref{fig:ICD_training} displays the convergence of voltage deviations and equity over 50 iterations for configuration 1. The equity value is also approaching zero (at a slower rate). We note that this plot shows the voltage deviations from the linearized system, while the deviations in Table~\ref{tab:v_compare} are reported from the power flow simulation software.

\begin{figure}[h!]
    \centering
    \includegraphics[width=.85\linewidth]{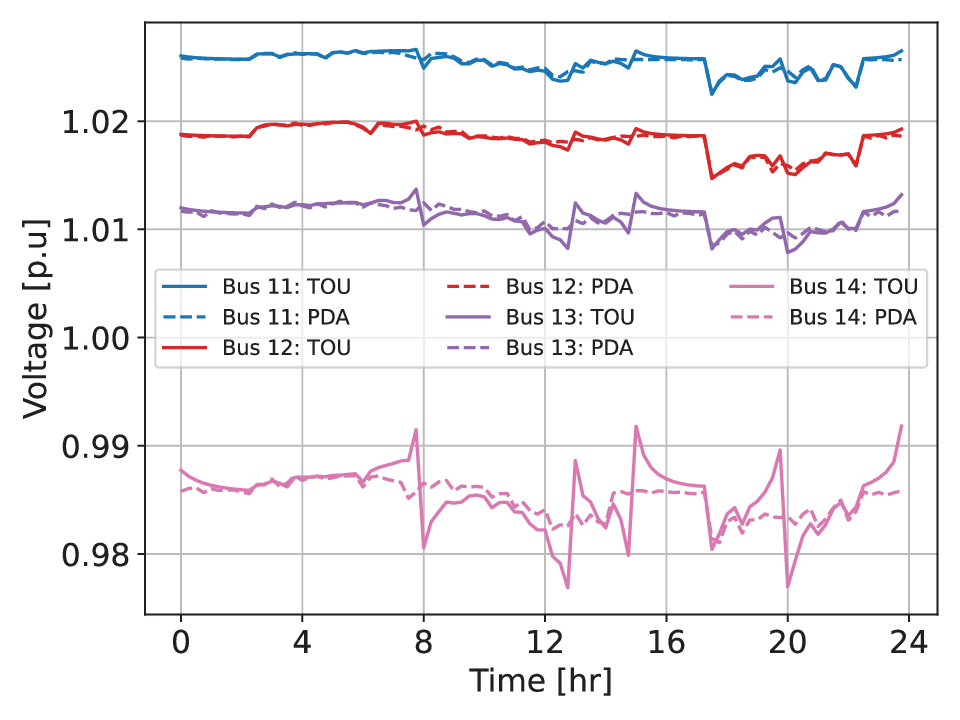}
    \caption{Voltages for buses 11-14 for SDID algorithm and configuration 1. Similar to ICD, the battery at bus 14 can smooth out the voltage profile. }
    \label{fig:sdid_powersim_1}
\end{figure}

\begin{figure}[h!]
    \centering
    \includegraphics[width=.85\linewidth]{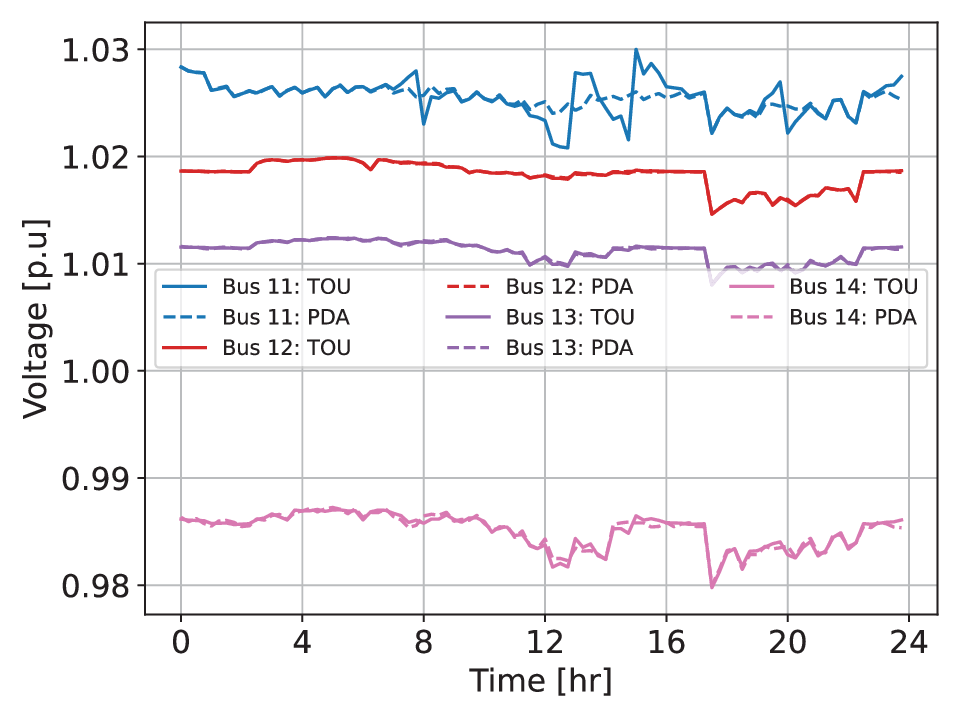}
    \caption{Voltages for buses 11-14 for SDID algorithm and configuration 2. The large battery at bus 11 can smooth out the voltage profile.}
    \label{fig:sdid_powersim_2}
\end{figure}

The plots in Figures~\ref{fig:sdid_powersim_1} and~\ref{fig:sdid_powersim_2} show that the response of the EVSP to prices addresses voltage deviations. In configuration 1, Bus 14 is assigned the biggest battery, leading to more effective voltage smoothing. Similarly, in configuration 2, bus 11 (which has the highest capacity battery) is able to effectively smooth out its voltage profile by responding to prices.

\subsection{ICD and SDID Comparison}
Table~\ref{tab:v_compare} provides a performance summary for both ICD and SDID. The reported voltage deviation represents the voltage fluctuations in per unit (p.u.) and is calculated by taking the Frobenius norm of the $\Delta \mathbf{V}$ matrix, which is the voltage deviation over the time horizon. The bold values represent the best performing algorithm for each metric. ICD and SDID yield the minimum voltage deviation in configuration 1 and 2, respectively, with both showing at least 45\% improvement over the baseline. ICD always outperforms SDID in terms of equity, reporting equity values very close to the global optimum of 0. ICD is computationally faster (0.4 seconds per iteration) than SDID (8 seconds per iteration). This is because SDID requires additional operations for obtaining the total gradient of the objective with respect to prices, while solving a similar EVSP profit maximization subproblem as ICD every iteration.
Additionally, ICD produces effective equitable prices after just one iteration, while SDID requires at least 50 iterations for similar convergence (see Fig.~\ref{fig:ICD_training}). In terms of robustness, ICD converges to the best solution in 1-4 iterations, while SDID convergence depends on price initialization, step size, and stepping method. Thus, ICD tends to be more robust with fewer hyperparameters to tune.


\begin{table}[h]
\setlength{\tabcolsep}{4pt}
\caption{\hspace{0.1em}Results summary}
\centering
\begin{tabular}{ |c|c|c|c|c|c|c|c|} 
\hline
& \multicolumn{3}{c}{Voltage Deviation (p.u)} & \multicolumn{2}{|c|}{Equity}\\
\hline
Config. & TOU (baseline) & ICD & SDID & ICD & SDID\\
\hline 
 1 & $0.0337$ & $\mathbf{0.0153}$ & $0.0187$ & $\mathbf{-1.27e^{-6}}$ & $-0.085$\\
\hline
 2 & $0.0249$ & $0.0197$ & $\mathbf{0.0169}$ & $\mathbf{-2.88e^{-6}}$ & $-0.823$\\
\hline
\end{tabular}
\label{tab:v_compare}
\end{table}

\section{Conclusion}\label{conclusion}
This paper presents a novel algorithm to set dynamic nodal electricity prices in a way that minimizes voltage violations and ensures equitable pricing for EVSPs in a distribution grid. We examine the performance of two proposed methods, ICD and SDID, by running power flow simulations on an IEEE 14-bus test system. Both methods outperform an existing PG\&E TOU rate. Compared to the SDID approach, the ICD approach is more computationally efficient, more robust to price initialization point, and has fewer hyperparameters.

In this work, EVSP agents are assumed to be honest; however, this is not guaranteed. Future work will integrate methods for addressing dishonest actors/market participants. The algorithms can also be applied to larger networks with other time-varying profiles. Specifically, the major energy ``producing'' devices in this work are limited to base generators and battery systems, but DERs such as solar PV and vehicle-to-grid (V2G) can be included in the setting. Deeper DER penetration may increase the variance of resource availability within the grid, increasing voltage fluctuations and exacerbating price inequity. Both the ICD and SDID algorithms can be used to examine this phenomenon.

Equity is vital to a just clean energy transition, and our proposed framework and algorithms provide a methodology for equitable electricity pricing without sacrificing power quality.

\bibliographystyle{IEEEtran}
\bibliography{IEEEfull}


\end{document}